# SMALL EIGENVALUES OF THE CONFORMAL LAPLACIAN

CHRISTIAN BÄR AND MATTIAS DAHL

ABSTRACT. We introduce a differential topological invariant for compact differentiable manifolds by counting the small eigenvalues of the Conformal Laplace operator. This invariant vanishes if and only if the manifold has a metric of positive scalar curvature. We show that the invariant does not increase under surgery of codimension at least three and we give lower and upper bounds in terms of the $\alpha$-genus.

## 1. INTRODUCTION

Throughout the paper let $M$ be a compact oriented differentiable manifold of dimension $n \geq 3$. Given a Riemannian metric $g$ on $M$ the Conformal Laplacian $L_g$ is defined as

$$L_g = \Delta_g + \frac{n-2}{4(n-1)} \cdot \text{Scal}_g$$

where $\Delta_g = d^*d$ is the Laplacian and $\text{Scal}_g$ is the scalar curvature of $g$. The operator $L_g$ is an elliptic differential operator of second order, self-adjoint in $L^2(M, \mathbb{R})$. Let $\mu_0(L_g) \leq \mu_1(L_g) \leq \mu_2(L_g) \leq \ldots$ be the spectrum of $L_g$, the eigenvalues being repeated according to their multiplicities. Let $f$ be a positive function on $M$. The Conformal Laplacian of the conformally related metric $\overline{g} = f^{\frac{4}{n-2}} g$ is given by

$$(1) \qquad L_{\overline{g}} u = f^{-\frac{n+2}{n-2}} L_g(fu).$$

Applying (1) to the function $u = 1$ gives the formula

$$(2) \qquad \text{Scal}_{\overline{g}} = \frac{4(n-1)}{n-2} f^{-\frac{n+2}{n-2}} L_g f$$

for the scalar curvature of $\overline{g}$.

We now introduce a differential topological invariant of a compact manifold by counting the number of small eigenvalues of the Conformal Laplacian.

**Definition 1.1.** Let $M$ be a compact differentiable manifold. The *$\kappa$-invariant* $\kappa(M)$ is defined to be the smallest integer $k$ such that for every $\varepsilon > 0$ there is a Riemannian metric $g_\varepsilon$ on $M$ for which

$$\begin{cases} \mu_k(L_{g_\varepsilon}) = 1, \\ |\mu_i(L_{g_\varepsilon})| < \varepsilon, \quad 0 \leq i < k. \end{cases}$$

2000 *Mathematics Subject Classification.* 53C27, 55N22, 57R65, 58J05, 58J50.

The first author was partially supported by the Research and Training Network "EDGE". The second author was supported by a Post-doctoral scholarship from the Swedish Foundation for International Cooperation in Research and Higher Education "STINT". Both authors have also been partially supported by the Research and Training Network "Geometric Analysis" funded by the European Commission.

If no such integer exists set $\kappa(M) := \infty$.

Heuristically, $\kappa(M)$ is the dimension of the "almost-kernel" of the Conformal Laplace operator.

By rescaling the metrics $g_\varepsilon$ accordingly one sees that $\kappa(M)$ is also the smallest integer $k$ such that for each constant $C > 0$ there exists a Riemannian metric $g_C$ for which

$$\begin{cases} \mu_k(L_{g_C}) > C, \\ |\mu_i(L_{g_C})| \leq 1, \quad 0 \leq i < k. \end{cases}$$

Hence $\kappa(M)$ tells us which is the first eigenvalue that can be made arbitrarily large for appropriate metrics while keeping the preceeding ones bounded.

If we made this definition using the Laplace operator acting on $p$-forms instead of the Conformal Laplacian, then by Hodge theory the resulting invariant would be nothing but the $p^{\text{th}}$ Betti number.

From the fact that the spectrum of $M_1 + M_2$ is the disjoint union of the spectra of $M_1$ and of $M_2$ it follows that

(3) $$\kappa(M_1 + M_2) = \kappa(M_1) + \kappa(M_2)$$

where we use sum notation to denote disjoint unions of manifolds. Also $\kappa(-M) = \kappa(M)$, where $-M$ denotes $M$ with reversed orientation. The next proposition concerns the relation between $\kappa(M)$ and scalar curvature.

**Proposition 1.2.** *Let $M$ be a compact differentiable manifold of dimension $n \geq 3$. Then*
  (1) *$\kappa(M) = 0$ if and only if there is a metric of positive scalar curvature on $M$.*
  (2) *If $M$ is connected and has a scalar-flat metric then $\kappa(M) \leq 1$.*

*Proof.* If $\kappa(M) = 0$ then there is a metric with $\mu_0 = 1$. The corresponding eigenfunction $f_0$ can chosen to be positive. From Equation (2) it follows that $\overline{g} = f_0^{\frac{4}{n-2}} g$ has positive scalar curvature. Conversely, if $g$ is a metric of positive scalar curvature on $M$ then $L_g > 0$ and we can rescale so that $\mu_0 = 1$. Hence $\kappa(M) = 0$.

For a scalar-flat metric on $M$ we have $L = \Delta \geq 0$ and the zero eigenspace consists of the constant functions. So if $M$ is connected we have $\mu_0 = 0$ and $\mu_1 > 0$. □

The following theorem controls the spectrum of $L_g$ under surgeries of codimension at least three. This will enable us to examine the behavior of $\kappa(M)$ under such surgeries.

**Theorem 3.1.** *Let $(M, g)$ be a closed Riemannian manifold. Let $\widetilde{M}$ be obtained from $M$ by surgery in codimension at least three. Then for each $k \in \mathbb{N}$ and for each $\varepsilon > 0$ there exists a Riemannian metric $\widetilde{g}$ on $\widetilde{M}$ such that the first $k + 1$ eigenvalues of the operators $L_g$ and $L_{\widetilde{g}}$ are $\varepsilon$-close, that is*

$$|\mu_j(L_g) - \mu_j(L_{\widetilde{g}})| < \varepsilon$$

*for $j = 0, \ldots, k$.*

As an immediate consequence we obtain

**Corollary 3.2.** *Let $M$ be a compact differentiable manifold of dimension $n \geq 3$. Suppose $\widetilde{M}$ is obtained from $M$ by surgery of codimension $\geq 3$. Then*
$$\kappa(\widetilde{M}) \leq \kappa(M).$$

Hence for any $\kappa_0 \in \mathbb{N}_0$ the property of having $\kappa \leq \kappa_0$ is preserved under surgery of codimension at least three. For $\kappa_0 = 0$ this means that the property of admitting a metric of positive scalar curvature is preserved under such surgeries. This is a famous by now classical result of Gromov and Lawson [10]. We do not give a new proof of this fact since we use the work of Gromov and Lawson when we prove Theorem 3.1.

As to the case $\kappa_0 = 1$ it is interesting to note that the property of allowing a scalar flat metric is not preserved under such surgeries. It follows that the converse of statement (2) in Proposition 1.2 does not hold. For example, the $n$-dimensional torus $T^n$ has a flat metric but no metric of positive scalar curvature [11]. Thus $\kappa(T^n) = 1$. Performing surgery in codimension at least three on $T^n$ yields a manifold $M^n$ not admitting metrics with positive or zero scalar curvature. Yet we have $\kappa(M^n) = 1$.

Also note that the condition $\kappa = 0$ is not preserved under surgery of codimension 2. Like any compact connected 3-manifold the 3-torus $T^3$ can be obtained from $S^3$ by a sequence of surgeries in codimension 2. But we have $\kappa(T^3) = 1 > \kappa(S^3) = 0$. This also shows that Theorem 3.1 cannot hold for surgeries in codimension less than three.

The $\kappa$-invariant measures how close $L$ can come to being a positive operator for some Riemannian metric on $M$. Since $L$ is positive if and only if $M$ allows a metric of positive scalar curvature one can also view $\kappa$ as a measure of how close one can get to having positive scalar curvature. Therefore it is not unreasonable to suspect that $\kappa$ is related to the $\hat{A}$ or $\alpha$-genus of $M$, the primary obstruction to allowing metrics of positive scalar curvature. We will see that this indeed is the case. On the one hand we have

**Theorem 2.4.** *Let $M$ be a compact spin manifold of dimension $n = 4m$. Then*
$$|\hat{A}(M)| \leq 2^{2m-1}\kappa(M).$$

As an application we obtain the following isoperimetric result.

**Corollary 2.5.** *Let $M$ be a compact spin manifold of dimension $n = 4m$ with $|\hat{A}(M)| > 2^{2m-1}$. Then there exists a constant $C = C(M)$ such that for each Riemannian metric with $|\operatorname{Scal}| \leq 1$ there exists a hypersurface $S \subset M$ dividing $M$ into two connected components $M_1$ and $M_2$ such that*
$$\operatorname{vol}_{n-1}(S) \leq C \cdot \min\{\operatorname{vol}_n(M_1), \operatorname{vol}_n(M_2)\}.$$

On the other hand we can bound $\kappa(M)$ from above in terms of the dimension and the $\alpha$-genus, at least for simply connected manifolds of dimension $n \geq 5$. First we make the following

**Observation.** *Let $M$ be a simply connected compact differentiable manifold of dimension $n \geq 5$. If $M$ is non-spin or if $n \equiv 3, 5, 6, 7 \mod 8$ then*
$$\kappa(M) = 0.$$

This comes from the fact that in these cases $M$ is well-known to carry a metric of positive scalar curvature, see [10], [17].

In dimensions $n \equiv 0 \mod 4$ the $\alpha$-genus of a spin manifold is integer-valued and it essentially coincides with the $\hat{A}$-genus. More precisely, if $n = 8l$ then $\alpha(M) = \hat{A}(M)$ while if $n = 8l + 4$ then $\alpha(M) = \frac{1}{2}\hat{A}(M)$.

**Theorem 4.4** *Let $M$ be a simply connected differentiable manifold of dimension $n \equiv 0 \mod 4$. Write $n = 8l$ or $n = 8l + 4$ with $l \geq 1$ and let $|\alpha(M)| = 4^l p + q$, $p \geq 0$, $0 \leq q < 4^l$. Then*

$$\kappa(M) \leq p + \min\{q, l\}.$$

As a special case we see that for spin manifolds as in the Theorem we have $\kappa = 1$ if $\alpha = 1$. In dimensions $n \equiv 1, 2 \mod 8$ this and the converse is true. In those dimensions we have $\alpha(M) \in \mathrm{KO}^{-n}(\mathrm{pt}) \cong \mathbb{Z}/2\mathbb{Z}$. By $|\alpha(M)| \in \mathbb{Z}$ we mean 0 if $\alpha(M)$ is trivial and 1 otherwise.

**Theorem 4.6.** *Let $M$ be a simply connected spin manifold of dimension $n = 8l + 1$ or $8l + 2$, $l \geq 1$. Then*

$$\kappa(M) = |\alpha(M)|.$$

This shows that $\kappa(M)$ can distinguish certain exotic spheres. In particular, $\kappa(M)$ is *not* invariant under homeomorphisms, only under diffeomorphisms.

Even though Theorem 2.4 shows that $\kappa(M)$ can become arbitrarily large it turns out that in a stable sense it takes only the values 0 and 1. More precisely, let $B$ be a compact simply connected 8-dimensional spin manifold with $\hat{A}(B) = 1$. Then $\alpha(M \times B) = \alpha(M)$ for all spin manifolds $M$.

**Theorem 5.1.** *Let $M$ be a simply connected spin manifold. Then*

$$\kappa(M \times B^p) \leq 1$$

*for all sufficiently large $p$.*

In the interesting work [15] the scalar curvature-related Yamabe invariant is studied using similar applications of the surgery and bordism results of [10], [17].

The paper is organized as follows. In Section 2 we estimate $\kappa(M)$ from below and prove Theorem 2.4. This is achieved by comparing $\kappa(M)$ to the dimension of the kernel of the Dirac operator using a spectral comparison principle of Gallot and Meyer. This is another manifestion of a deep relationship between the spectrum of the Conformal Laplacian and the Dirac spectrum which was first observed by Hijazi who compares the lowest eigenvalues in [12]. Another important ingredient is a refined Kato inequality for harmonic spinors. In Section 3 we study the behavior of the spectrum of the Conformal Laplacian and prove Theorem 3.1. This together with standard results from bordism theory is used in Section 4 to derive upper bounds on $\kappa(M)$. Manifolds with special holonomy occur as important building blocks. In the final section we study the stable limit and prove Theorem 5.1. Moreover, we discuss the behavior of $\kappa(M)$ when $M$ is replaced by a finite covering. It turns out that not much can be said in general, $\kappa$ can decrease, increase or remain unaltered.

The authors wish to thank Stephan Stolz and the referee for insightful comments.

## 2. Lower bound for $\kappa$

In this section we are going to find a lower bound on $\kappa$ in terms of the $\hat{A}$-genus. This will follow from a spectral comparison result relating the kernel of the Dirac operator to the spectrum of the Conformal Laplacian, in the spirit of Gallot and Meyer [8]. The two main technical points of the spectral comparison theorem are the "Hilbertian lemma" of [8] and a refined Kato inequality for the Dirac operator, which we now recall.

Let $E$ be a Riemannian vector bundle over $M$ of real fiber dimension $l$. Let $H$ be a subspace of $L^2(M, E)$ of (real) dimension $h$ and let $K$ be a subspace of $L^2(M, \mathbb{R})$ of dimension $k$. Let $\pi_K$ be the orthogonal projection onto $K$.

**Theorem 2.1.** [8] *Suppose $h \geq k(l+1)$. Then there exists $\varphi \in H$ such that*

$$\|\pi_K(|\varphi|)\|^2 < \left(1 - \frac{1}{8(l+1)^2}\right)\|\varphi\|^2. \tag{4}$$

The classical Kato inequality states that $|d|\varphi||^2 \leq |\nabla\varphi|^2$ wherever $\varphi \neq 0$. If $\varphi$ is a harmonic spinor field this can be improved as follows. Denote the (complex) spinor bundle of a Riemannian spin manifold $M$ by $\Sigma M$.

**Proposition 2.2.** *Suppose $\varphi$ is a harmonic spinor on an $n$-dimensional Riemannian spin manifold, i. e. $D\varphi = 0$ where $D$ is the Dirac operator. Then*

$$|d|\varphi||^2 \leq \frac{n-1}{n}|\nabla\varphi|^2$$

*at all points where $\varphi$ is non-zero.*

*Proof.* This proposition is a special case of the general work on refined Kato inequalities in [3], [4]. For the convenience of the reader we give a simple direct proof. Fix a point $p \in M$ at which $\varphi(p) \neq 0$ so that $|\varphi|$ is differentiable at $p$. We define $\pi \in \mathrm{End}(T_pM \otimes \Sigma_pM)$ by

$$\pi(X \otimes \psi) := -\frac{1}{n}\sum_{j=1}^n e_j \otimes e_j \cdot X \cdot \psi,$$

where $e_1, \ldots, e_n$ is an orthonormal basis of $T_pM$ and $\cdot$ denotes Clifford multiplication. Here we use the real tensor product of the real vector space $T_pM$ and the complex vector space $\Sigma_pM$ to obtain a *complex* vector space $T_pM \otimes \Sigma_pM$. Equivalently, we could complexify $T_pM$ and then use the complex tensor product. Moreover, $T_pM \otimes \Sigma_pM$ inherits a Hermitian scalar product.

It is readily checked that $\pi$ is independent of the choice of orthonormal basis and that $\pi$ is an orthogonal projection. Denote the complementary projection $1 - \pi$ by $\pi'$. We compute

$$
\begin{aligned}
\langle \pi'(X \otimes \psi), \pi'(X \otimes \psi) \rangle &= \langle X \otimes \psi, X \otimes \psi \rangle \\
&\quad - 2\operatorname{Re}\left\langle X \otimes \psi, -\frac{1}{n}\sum_j e_j \otimes e_j \cdot X \cdot \psi \right\rangle \\
&\quad + \left\langle -\frac{1}{n}\sum_j e_j \otimes e_j \cdot X \cdot \psi, -\frac{1}{n}\sum_k e_k \otimes e_k \cdot X \cdot \psi \right\rangle \\
&= |X|^2|\psi|^2 + \frac{2}{n}\sum_j \operatorname{Re}\langle X, e_j\rangle \langle \psi, e_j \cdot X \cdot \psi\rangle \\
&\quad + \frac{1}{n^2}\sum_j \langle e_j \cdot X \cdot \psi, e_j \cdot X \cdot \psi \rangle \\
&= |X|^2|\psi|^2 - \frac{2}{n}|X|^2|\psi|^2 + \frac{1}{n}|X|^2|\psi|^2 \\
&= \frac{n-1}{n}|X|^2|\psi|^2.
\end{aligned}
$$

Now we return to the harmonic spinor $\varphi$. Note that

$$
\pi(\nabla\varphi) = \pi\left(\sum_k e_k \otimes \nabla_{e_k}\varphi\right) = -\frac{1}{n}\sum_{j,k} e_j \otimes e_j e_k \nabla_{e_k}\varphi = -\frac{1}{n}\sum_j e_j \otimes e_j D\varphi = 0,
$$

so $\nabla\varphi = \pi'(\nabla\varphi)$. Choose a unit vector $X \in T_pM$ such that $|d|\varphi|| = \langle X, \operatorname{grad}|\varphi|\rangle$. Then we have

$$
\begin{aligned}
|\varphi||d|\varphi|| &= \left\langle X, \frac{1}{2}\operatorname{grad}\langle \varphi, \varphi\rangle \right\rangle \\
&= \operatorname{Re}\langle \nabla_X\varphi, \varphi\rangle \\
&\leq |\langle \nabla_X\varphi, \varphi\rangle| \\
&= |\langle \nabla\varphi, X \otimes \varphi\rangle| \\
&= |\langle \pi'(\nabla\varphi), X \otimes \varphi\rangle| \\
&= |\langle \nabla\varphi, \pi'(X \otimes \varphi)\rangle| \\
&\leq |\nabla\varphi||\pi'(X \otimes \varphi)| \\
&\leq |\nabla\varphi|\left(\frac{n-1}{n}\right)^{1/2}|X||\varphi| \\
&= \left(\frac{n-1}{n}\right)^{1/2}|\varphi||\nabla\varphi|.
\end{aligned}
$$

$\square$

For $\varepsilon > 0$ define the smooth approximation $|\cdot|_\varepsilon$ of the norm $|\cdot|$ by $|\varphi|_\varepsilon := \left(|\varphi|^2 + \varepsilon^2\right)^{1/2}$. If $\varphi$ is harmonic it follows from the refined Kato inequality that

(5)    $$|d|\varphi|_\varepsilon|^2 = \frac{|\varphi|^2}{|\varphi|_\varepsilon^2}|d|\varphi||^2 \leq \frac{n-1}{n}|\nabla\varphi|^2$$

at points where $\varphi \neq 0$. Since $\{x \in M \mid \varphi(x) \neq 0\}$ is dense in $M$ for harmonic $\varphi$ we conclude that (5) holds on all of $M$.

Now let $M$ be a compact Riemannian spin manifold and let $SM$ be a sum of copies of the spinor bundle $\Sigma M$, or in even dimension a sum of copies of $\Sigma^+ M$ and $\Sigma^- M$, the bundles of positive and of negative half-spinors. Let $l$ be the real rank of $SM$ and let $h$ be the real dimension of the kernel of $D^2$ acting on sections of $SM$. Let $C(n,l) = \frac{8(l+1)^2(n-1)^2}{n(n-2)} - 1$.

**Theorem 2.3.** *Suppose that for some $k \geq 1$ we have*
$$\mu_k(L_g) \geq -C(n,l)\mu_0(L_g).$$
*Then*
$$h < k(l+1).$$

*Proof.* Let $K$ be the subspace of $L^2(M,\mathbb{R})$ spanned by the eigenfunctions of $L_g$ corresponding to the eigenvalues $\mu_0, \mu_1, \ldots, \mu_{k-1}$ and let $H \subset L^2(M, SM)$ be the space of harmonic spinors, so that $\dim_\mathbb{R} H = h$.

Let $f$ be a positive eigenfunction for $\mu_0$ and suppose $h \geq k(l+1)$. Theorem 2.1 then gives us a spinor in $f^{-\frac{1}{n-2}} H$ satisfying (4), that is

(6) $$\|\pi_K(|f^{-\frac{1}{n-2}}\varphi|)\|^2 < \left(1 - \frac{1}{8(l+1)^2}\right)\|f^{-\frac{1}{n-2}}\varphi\|^2$$

where $\varphi \in H$.

Define the conformally related metric $\overline{g} = f^{\frac{4}{n-2}}g$. Then $\overline{\varphi} = f^{-\frac{n-1}{n-2}}\varphi$ satisfies $\overline{D}\overline{\varphi} = 0$ where $\overline{D}$ is the Dirac operator defined using $\overline{g}$, see for instance [2]. From the Schrödinger-Lichnerowicz formula and (5) we have

$$\begin{aligned}
0 &= \frac{n-1}{n}\|\overline{D}\overline{\varphi}\|_{\overline{g}}^2 \\
&= \frac{n-1}{n}\|\overline{\nabla}\overline{\varphi}\|_{\overline{g}}^2 + \frac{n-1}{4n}\left(\mathrm{Scal}_{\overline{g}}\,\overline{\varphi}, \overline{\varphi}\right)_{\overline{g}} \\
&\geq \left(\overline{\Delta}|\overline{\varphi}|_\varepsilon, |\overline{\varphi}|_\varepsilon\right)_{\overline{g}} + \frac{n-1}{4n}\left(\mathrm{Scal}_{\overline{g}}\,\overline{\varphi}, \overline{\varphi}\right)_{\overline{g}} \\
&= (L_{\overline{g}}|\overline{\varphi}|_\varepsilon, |\overline{\varphi}|_\varepsilon)_{\overline{g}} - \frac{n-2}{4(n-1)}\left(\mathrm{Scal}_{\overline{g}}\,|\overline{\varphi}|_\varepsilon, |\overline{\varphi}|_\varepsilon\right)_{\overline{g}} + \frac{n-1}{4n}\left(\mathrm{Scal}_{\overline{g}}\,|\overline{\varphi}|, |\overline{\varphi}|\right)_{\overline{g}}.
\end{aligned}$$

Using (1) we compute for the first term

$$\begin{aligned}
(L_{\overline{g}}|\overline{\varphi}|_\varepsilon, |\overline{\varphi}|_\varepsilon)_{\overline{g}} &= \int_M f^{-\frac{n+2}{n-2}} L_g(f|\overline{\varphi}|_\varepsilon)|\overline{\varphi}|_\varepsilon f^{\frac{2n}{n-2}}\, dV_g \\
&= (L_g(f|\overline{\varphi}|_\varepsilon), f|\overline{\varphi}|_\varepsilon)_g \\
&= (L_g \pi_K(f|\overline{\varphi}|_\varepsilon), \pi_K(f|\overline{\varphi}|_\varepsilon))_g + (L_g \pi_{K^\perp}(f|\overline{\varphi}|_\varepsilon), \pi_{K^\perp}(f|\overline{\varphi}|_\varepsilon))_g \\
&\geq \mu_0 \|\pi_K(f|\overline{\varphi}|_\varepsilon)\|^2 + \mu_k\left(\|f|\overline{\varphi}|_\varepsilon\|^2 - \|\pi_K(f|\overline{\varphi}|_\varepsilon)\|^2\right) \\
&= (\mu_0 - \mu_k)\|\pi_K(f|\overline{\varphi}|_\varepsilon)\|^2 + \mu_k \|f|\overline{\varphi}|_\varepsilon\|^2
\end{aligned}$$

where $\pi_K$ and $\pi_{K^\perp}$ are the orthogonal projections onto $K$ and $K^\perp$. Letting $\varepsilon$ go to zero we get

(7) $$0 \geq (\mu_0 - \mu_k)\|\pi_K(f|\overline{\varphi}|)\|^2 + \mu_k\|f|\overline{\varphi}|\|^2 + \frac{1}{4n(n-1)}\left(\mathrm{Scal}_{\overline{g}}\,|\overline{\varphi}|, |\overline{\varphi}|\right)_{\overline{g}}.$$

Using (6) we can estimate the first term of (7) as

$$\begin{aligned}(\mu_0 - \mu_k)\|\pi_K(f|\overline{\varphi}|)\|^2 &= (\mu_0 - \mu_k)\|\pi_K(|f^{-\frac{1}{n-2}}\varphi|)\|^2 \\ &> (\mu_0 - \mu_k)\left(1 - \frac{1}{8(l+1)^2}\right)\|f^{-\frac{1}{n-2}}\varphi\|^2.\end{aligned}$$

From (2) we get for the third term of (7)

$$\begin{aligned}(\mathrm{Scal}_{\overline{g}}|\overline{\varphi}|, |\overline{\varphi}|)_{\overline{g}} &= \int_M \mathrm{Scal}_{\overline{g}} f^{-2\frac{n-1}{n-2}}|\varphi|^2 f^{\frac{2n}{n-2}} dV_g \\ &= \int_M \mathrm{Scal}_{\overline{g}} f^{\frac{4}{n-2}}|f^{-\frac{1}{n-2}}\varphi|^2 dV_g \\ &= \frac{4(n-1)}{n-2}\int_M f^{-1}L_g f|f^{-\frac{1}{n-2}}\varphi|^2 dV_g \\ &= \frac{4(n-1)}{n-2}\mu_0\|f^{-\frac{1}{n-2}}\varphi\|^2.\end{aligned}$$

Together we have

$$\begin{aligned}0 &> \left((\mu_0 - \mu_k)\left(1 - \frac{1}{8(l+1)^2}\right) + \mu_k + \frac{1}{n(n-2)}\mu_0\right)\|f^{-\frac{1}{n-2}}\varphi\|^2 \\ &= \frac{1}{8(l+1)^2}\left(\mu_k + \left(\frac{8(l+1)^2(n-1)^2}{n(n-2)} - 1\right)\mu_0\right)\|f^{-\frac{1}{n-2}}\varphi\|^2 \\ &= \frac{1}{8(l+1)^2}\left(\mu_k + C(n,l)\mu_0\right)\|f^{-\frac{1}{n-2}}\varphi\|^2\end{aligned}$$

so $\mu_k < -C(n,l)\mu_0$, which contradicts the assumption in the Theorem. We conclude that $h < k(l+1)$. □

Using Theorem 2.3 we now prove a lower bound for $\kappa(M)$.

**Theorem 2.4.** *Let $M$ be a compact spin manifold of dimension $n = 4m$. Then*

$$|\hat{A}(M)| \leq 2^{2m-1}\kappa(M).$$

*Proof.* If $\kappa(M) = 0$ then $M$ has a metric of positive scalar curvature and $\hat{A}(M) = 0$, so we assume $\kappa(M) > 0$. Also we assume that $\hat{A}(M) \geq 0$; if this is not the case then change the orientation of $M$.

Let $p$ be an integer such that $2p > \kappa(M)$. Choose a Riemannian metric $g$ on $M$ for which $|\mu_0(L_g)| \leq 1$ and $\mu_{\kappa(M)}(L_g) > C(n,l)$ where $l$ is the real rank of $SM = p\Sigma^+M$. Hence

$$\mu_{\kappa(M)}(L_g) > -C(n,l)\mu_0(L_g).$$

We have $\hat{A}(M) = h_+ - h_- \leq h_+$ where $h_\pm$ is the complex dimension of the kernel of $D^2$ acting on sections of $\Sigma^\pm M$. Consider $D^2$ acting on sections of $SM$. The real dimension of the kernel is $p \cdot 2 \cdot h_+$ and the real rank of $SM$ is $p \cdot 2 \cdot 2^{[n/2]-1} = p \cdot 2^{2m}$. From Theorem 2.3 we conclude that

$$p \cdot 2 \cdot h_+ < \kappa(M)\left(p \cdot 2^{2m} + 1\right)$$

so

$$|\hat{A}(M)| = \hat{A}(M) \leq h_+ < 2^{2m-1}\kappa(M) + \frac{\kappa(M)}{2p}.$$

The result follows since $\frac{\kappa(M)}{2p} < 1$ and the other terms are integers. □

Since manifolds admitting Riemannian metrics with non-negative scalar curvature satisfy $\kappa(M) \leq 1$ we get $|\hat{A}(M)| \leq 2^{2m-1}$ for such manifolds. In case $M$ has finite fundamental group $\pi_1(M)$ Futaki [7, Cor. 2] shows the stronger estimate $|\pi_1(M)| \cdot |\hat{A}(M)| \leq 2^m$ for manifolds admitting non-negative scalar curvature. If $M$ has infinite fundamental group then Mathai [14, Thm. 1.3] proves that $\hat{A}(M) = 0$.

As an application of Theorem 2.4 we obtain the following result saying that under the topological condition $|\hat{A}(M)| > 2^{n/2-1}$ and under the normalizing condition $|\operatorname{Scal}| \leq 1$ the manifold has a "neck" of uniformly bounded size.

**Corollary 2.5.** *Let $M$ be a compact spin manifold of dimension $n = 4m$ with $|\hat{A}(M)| > 2^{2m-1}$. Then there exists a constant $C = C(M)$ such that for each Riemannian metric with $|\operatorname{Scal}| \leq 1$ there exists a hypersurface $S \subset M$ dividing $M$ into two connected components $M_1$ and $M_2$ such that*

$$\operatorname{vol}_{n-1}(S) \leq C \cdot \min\{\operatorname{vol}_n(M_1), \operatorname{vol}_n(M_2)\}.$$

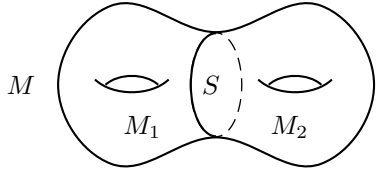

FIG. 1

*Proof.* From $|\hat{A}(M)| > 2^{2m-1}$ we conclude by Theorem 2.4 that $\kappa(M) > 1$. Thus there exists a constant $C_1$ such that for each metric $g$ we have $|\mu_0(L_g)| > \frac{4(n-1)}{n-2}$ or $\mu_1(L_g) \leq C_1$. For a metric $g$ with $|\operatorname{Scal}_g| \leq 1$ we obviously have $|\mu_0(L_g)| \leq \frac{4(n-1)}{n-2}$ and hence we must have $\mu_1(L_g) \leq C_1$. This implies for the first positive Laplace eigenvalue $\mu_1(\Delta_g) \leq C_1 + \frac{4(n-1)}{n-2}$. Cheeger's inequality [5] says that $\frac{\mathfrak{h}^2}{4} \leq \mu_1(\Delta_g)$ where

$$\mathfrak{h} = \inf_S \frac{\operatorname{vol}_{n-1}(S)}{\min\{\operatorname{vol}_n(M_1), \operatorname{vol}_n(M_2)\}}.$$

The corollary follows. □

Note that this statement is not true for manifolds admitting scalar flat metrics since one can then rescale the metric without violating $|\operatorname{Scal}| \leq 1$ and make Cheeger's constant $\mathfrak{h}$ arbitrarily large. Hence an assumption like $|\hat{A}(M)| > 2^{2m-1}$ (or rather $\kappa(M) > 1$) is necessary.

## 3. SURGERY AND THE CONFORMAL LAPLACIAN

The aim of this section is to study the behavior of the spectrum of the operator

$$L_g = \Delta_g + c \cdot \operatorname{Scal}_g$$

under surgery on the underlying manifold. Here $c$ is a fixed positive constant which could e. g. be $c = \frac{n-2}{4(n-1)}$. Since the precise value of $c$ is irrelevant for the results of this section, in particular the conformal behavior of the operator plays no role, we will work in this slightly larger generality, still denoting the resulting operator by $L_g$. More specifically, we will show

**Theorem 3.1.** *Let $(M, g)$ be a closed Riemannian manifold. Let $\widetilde{M}$ be obtained from $M$ by surgery in codimension at least three. Then for each $k \in \mathbb{N}$ and for each $\varepsilon > 0$ there exists a Riemannian metric $\widetilde{g}$ on $\widetilde{M}$ such that the first $k + 1$ eigenvalues of the operators $L_g = \Delta_g + c \cdot \mathrm{Scal}_g$ and $L_{\widetilde{g}} = \Delta_{\widetilde{g}} + c \cdot \mathrm{Scal}_{\widetilde{g}}$ are $\varepsilon$-close, that is*

$$|\mu_j(L_g) - \mu_j(L_{\widetilde{g}})| < \varepsilon$$

*for $j = 0, \ldots, k$.*

As an immediate consequence of this theorem we see that $\kappa$ does not increase under surgeries of codimension at least three.

**Corollary 3.2.** *Let $M$ be a compact differentiable manifold of dimension $n \geq 3$. Suppose $\widetilde{M}$ is obtained from $M$ by surgery of codimension $\geq 3$. Then*

$$\kappa(\widetilde{M}) \leq \kappa(M).$$

If $M$ can be recovered from $\widetilde{M}$ by a surgery of codimension $\geq 3$ then $\kappa(\widetilde{M}) = \kappa(M)$. This can be done e. g. if $\widetilde{M}$ can be obtained from $M$ by surgery of codimension $k$, $3 \leq k \leq n-2$ where $n \geq 5$, or if $\widetilde{M} = M \# (S^k \times S^{n-k})$ where $n \geq 5$ and $0 \leq k \leq n$.

**Remark 3.3.** Corollary 3.2 together with Theorem 2.4 makes it possible to compute $\kappa(M)$ exactly for some manifolds $M$. Let K3 be the 4-dimensional K3-surface, this is connected and has a scalar-flat metric but no metric of positive scalar curvature, so $\kappa(\mathrm{K3}) = 1$. Let

$$M = \mathrm{K3} \# \ldots \# \mathrm{K3}, \quad k \text{ summands}.$$

Then Corollary 3.2 tells us that

$$\kappa(M) \leq \kappa(\mathrm{K3} + \cdots + \mathrm{K3}) = k \cdot \kappa(\mathrm{K3}) = k.$$

On the other hand Theorem 2.4 gives us

$$2\kappa(M) \geq \hat{A}(M) = k \cdot \hat{A}(\mathrm{K3}) = 2k$$

so that

$$\kappa(M) = k.$$

Note that Theorem 3.1 would not be very meaningful for the Laplace-Beltrami operator, i. e. for $c = 0$. By a result of Colin de Verdière [6] it is possible to prescribe any finite number of Laplace eigenvalues.

The proof of Theorem 3.1 requires some preparation. First we show that despite the fact that second derivatives of the metric enter into the scalar curvature the eigenvalues of the operator $L_g$ depend continuously on the metric $g$ with respect to the $C^1$-topology.

**Lemma 3.4.** *Let $(M, g)$ be a closed Riemannian manifold. For each $k \in \mathbb{N}$ and for each $\varepsilon > 0$ there exists a $C^1$-open neighborhood of $g$ in the space of smooth Riemannian metrics such that for each $g'$ in this neighborhood we have*

$$|\mu_j(L_g) - \mu_j(L_{g'})| < \varepsilon$$

*for $j = 0, \ldots, k$.*

*Proof.* The quadratic form corresponding to the operator $L_g$ is given by
$$u \mapsto (L_g u, u)_{L^2(M,g)} = (du, du)_{L^2(M,g)} + c(\mathrm{Scal}_g\, u, u)_{L^2(M,g)}$$
where $(\cdot, \cdot)_{L^2(M,g)}$ denotes the $L^2$-scalar product with respect to $g$. No derivatives of the metric enter into $(du, du)_{L^2(M,g)}$. Choose a partition of unity $\chi_1, \ldots, \chi_r$ each $\chi_j$ having its support in a coordinate chart. Then
$$(L_g u, u)_{L^2(M,g)} = (du, du)_{L^2(M,g)} + c \sum_{j=1}^{r} (\mathrm{Scal}_g\, u, \chi_j u)_{L^2(M,g)}.$$

In local coordinates
$$\mathrm{Scal}_g = \frac{1}{2} \sum_{\alpha\beta\gamma\delta} g^{\alpha\gamma} g^{\beta\delta} \left( \frac{\partial^2 g_{\alpha\delta}}{\partial x^\beta \partial x^\gamma} + \frac{\partial^2 g_{\beta\gamma}}{\partial x^\alpha \partial x^\delta} - \frac{\partial^2 g_{\alpha\gamma}}{\partial x^\beta \partial x^\delta} - \frac{\partial^2 g_{\beta\delta}}{\partial x^\alpha \partial x^\gamma} \right) + A(g, dg)$$

where $A(g, dg)$ is an algebraic expression in the $g_{\alpha\beta}$ and their first derivatives. The crucial point here is that the second derivatives enter *linearly* into the expression for $\mathrm{Scal}_g$. This allows us to reduce the degree of derivatives of $g$ by a partial integration:

$$
\begin{aligned}
(\mathrm{Scal}_g\, u, \chi_j u) &= \int \chi_j \mathrm{Scal}_g\, u^2 \sqrt{\det(g_{\mu\gamma})} dx^1 \cdots dx^n \\
&= \int \frac{1}{2} \sum_{\alpha\beta\gamma\delta} \Big( -\frac{\partial g_{\alpha\delta}}{\partial x^\gamma} \frac{\partial}{\partial x^\beta} - \frac{\partial g_{\beta\gamma}}{\partial x^\delta} \frac{\partial}{\partial x^\alpha} + \frac{\partial g_{\alpha\gamma}}{\partial x^\delta} \frac{\partial}{\partial x^\beta} + \frac{\partial g_{\beta\delta}}{\partial x^\gamma} \frac{\partial}{\partial x^\alpha} \\
&\quad + A(g, dg) \Big) \Big( g^{\alpha\gamma} g^{\beta\delta} \chi_j u^2 \sqrt{\det(g_{\mu\gamma})} \Big) dx^1 \cdots dx^n.
\end{aligned}
$$

Since $\chi_j$ has compact support in the coordinate chart there are no boundary terms. Thus $(L_g u, u)_{L^2(M,g)}$ is the integral of a quadratic expression in $u$ and its first derivatives with coefficients being algebraic expressions in $g$ and its first derivatives. Hence $(L_g u, u)_{L^2(M,g)} - (L_{g'} u, u)_{L^2(M,g')}$ is the integral of a quadratic expression in $u$ and its first derivatives with small coefficients if $g'$ is $C^1$-close to $g$. It follows that
$$\left| (L_g u, u)_{L^2(M,g)} - (L_{g'} u, u)_{L^2(M,g')} \right| \leq \varepsilon(g') \left( (L_g + k)u, u \right)_{L^2(M,g)}$$
for all $u \in C^\infty(M)$ where $\varepsilon(g') \to 0$ as $g' \xrightarrow{C^1} g$ and $k$ is a constant sufficiently large as to make the operator $L_g + k$ positive, for example $k = |\min \mathrm{Scal}_g| + 1$. The lemma now follows from the variational characterization of eigenvalues
$$\mu_j(L_g) = \inf_{\substack{V \subset C^\infty(M) \\ \dim V = j+1}} \sup_{\substack{u \in V \\ u \neq 0}} \frac{(L_g u, u)_{L^2(M,g)}}{\|u\|_{L^2(M,g)}}$$
and similarly for $\mu_j(L_{g'})$. $\square$

**Remark 3.5.** The proof of the lemma also yields a simultaneous eigenvalue comparison for all eigenvalues, not just for the $k+1$ first ones. Namely, for each $\varepsilon > 0$ there is a $C^1$-neighborhood of $g$ such that for all $g'$ in this neighborhood
$$(1-\varepsilon)\mu_j(L_g) - \varepsilon \leq \mu_j(L_{g'}) \leq (1+\varepsilon)\mu_j(L_g) + \varepsilon$$
holds for all $j \geq 0$. But we will not need this here.

Next we derive information about the distribution of the $L^2$-norm of eigenfunctions. It turns out that only a little bit of the $L^2$-norm of eigenfunctions corresponding to low eigenvalues is contained in a region of large scalar curvature.

**Lemma 3.6.** *Let $M$ be a closed Riemannian manifold. Let $\Lambda, S_0, S_1 \in \mathbb{R}$, $S_0 < S_1$. Suppose the scalar curvature of $M$ satisfies $\mathrm{Scal} \geq S_0$. Put*

$$M_+ := \{x \in M \mid \mathrm{Scal}(x) \geq S_1\}.$$

*For any smooth function $u$ satisfying*

$$(Lu, u)_{L^2(M)} \leq \Lambda^2 \|u\|^2_{L^2(M)}$$

*the following inequality holds:*

$$\int_{M_+} |u|^2 dV \leq \frac{\Lambda^2/c - S_0}{S_1 - S_0} \int_M |u|^2 dV.$$

*Proof.* We compute

$$\begin{aligned} 0 &\leq \int_M |du|^2 dV = \int_M u \Delta u\, dV \\ &= \int_M uLu\, dV - c \int_M \mathrm{Scal}\, |u|^2 dV \\ &\leq \Lambda^2 \int_M |u|^2 dV - cS_0 \int_{M \setminus M_+} |u|^2 dV - cS_1 \int_{M_+} |u|^2 dV. \end{aligned}$$

Hence

$$(cS_1 - \Lambda^2) \int_{M_+} |u|^2 dV \leq (\Lambda^2 - cS_0) \int_{M \setminus M_+} |u|^2 dV$$

and therefore

$$(cS_1 - cS_0) \int_{M_+} |u|^2 dV \leq (\Lambda^2 - cS_0) \int_M |u|^2 dV.$$

$\square$

We will also need some finer control on the distribution of the $L^2$-norm on annular regions. For a compact submanifold $N \subset M$ and for $0 \leq R_1 < R_2$ define the "annular region"

$$A_N(R_1, R_2) := \{x \in M \mid R_1 \leq \mathrm{dist}(x, N) \leq R_2\}$$

and the distance sphere

$$S_N(R) := \{x \in M \mid \mathrm{dist}(x, N) = R\}.$$

**Lemma 3.7.** *Let $M$ be a Riemannian manifold and let $N \subset M$ be a compact submanifold of codimension at least 3. Then there exists $0 < R < 1$ such that for any $0 < r \leq R^{11}/2$ and any smooth function $u : A_N(r, (2r)^{1/11}) \to \mathbb{R}$ the following estimate holds*

$$\frac{\|u\|^2_{L^2(A_N(r,2r))}}{\|u\|^2_{L^2(A_N(r,(2r)^{1/11}))}} \leq 10\, r^{5/2}$$

*provided*

$$\int_{S_N(\rho)} u \partial_\nu u\, dA \geq 0$$

*holds for all $\rho \in [r, (2r)^{1/11}]$. Here $\nu$ denotes the unit normal vector field of $S_N(\rho)$ pointing away from $N$.*

The proof was given for spinors in [1], Lemma 2.4, and it carries over without changes to sections of an arbitrary Riemannian or Hermitean vector bundle equipped with a metric connection. In particular, it holds for functions.

*Proof of Theorem 3.1.* We introduce the eigenvalue counting function $\mathcal{N}_g : \mathbb{R} \to \mathbb{N}_0$ where $\mathcal{N}_g(\lambda)$ is the total number of eigenvalues $\leq \lambda$ of $L_g$ counted with multiplicity. We will show that given $\varepsilon > 0$ and $\Lambda \in \mathbb{R}$ there exists a metric $\widetilde{g}$ on $\widetilde{M}$ such that

$$\mathcal{N}_g(\lambda - \varepsilon) \leq \mathcal{N}_{\widetilde{g}}(\lambda) \leq \mathcal{N}_g(\lambda + \varepsilon)$$

for all $\lambda \leq \Lambda$. The theorem then follows easily.

Let $E_g(\lambda) \subset L^2(M)$ denote the direct sum of the $L_g$-eigenspaces for all eigenvalues $\leq \lambda$. Hence $\dim E_g(\lambda) = \mathcal{N}_g(\lambda)$. Let $N \subset M$ be the sphere along which surgery will be performed. Let $k$ be the codimension of $N$ in $M$. By assumption $k \geq 3$. Denote the distance tube of radius $r$ about $N$ by $U_N(r)$, that is,

$$U_N(r) := \{x \in M \mid \operatorname{dist}(x, N) < r\}.$$

We will first show that given $\varepsilon > 0$ and $\Lambda \in \mathbb{R}$ there exists $R > 0$ such that

$$\mathcal{N}_g(\lambda - \varepsilon) \leq \mathcal{N}_{\widetilde{g}}(\lambda)$$

for all $\lambda \leq \Lambda$ provided $(\widetilde{M}, \widetilde{g})$ contains an isometric copy of $M \setminus U_N(r)$ for some $0 < r \leq R$.

For $r > 0$ let $\chi_r : M \to \mathbb{R}$ be a smooth cut–off function such that

- $0 \leq \chi_r \leq 1$ on $M$,
- $\chi_r \equiv 0$ on $U_N(r)$,
- $\chi_r \equiv 1$ on $M \setminus U_N(2r)$,
- $|d\chi_r| \leq \frac{2}{r}$ on $M$.

Since $E_g(\Lambda)$ is finite dimensional there exists a constant $C_1 > 0$ such that

$$\|u\|_{L^\infty(M)} \leq C_1 \|u\|_{L^2(M)}$$

and

$$\|L_g u\|_{L^\infty(M)} \leq C_1 \|u\|_{L^2(M)}$$

for all $u \in E_g(\Lambda)$. Moreover, there is a constant $C_2 > 0$ such that

$$\operatorname{vol}(U_N(2r)) \leq C_2 r^k.$$

Now for $\lambda \leq \Lambda$ and $u \in E_g(\lambda - \varepsilon)$ the function $\chi_r u$ has its support in $M \setminus U_N(r)$ and can hence also be regarded as a function on $\widetilde{M}$. We plug it into the Rayleigh quotient of $L_{\widetilde{g}}$:

$$\begin{aligned}
(L_{\widetilde{g}}(\chi_r u), \chi_r u)_{L^2(\widetilde{M}, \widetilde{g})} &= (L_g(\chi_r u), \chi_r u)_{L^2(M, g)} \\
&= \int_M \left( \langle d(\chi_r u), d(\chi_r u) \rangle + c \operatorname{Scal}_g \chi_r^2 u^2 \right) dV \\
&= \int_M \left( \langle d(\chi_r u), d(\chi_r u) \rangle + \chi_r^2 u (L_g - \Delta) u \right) dV \\
&= \int_M \left( \langle d(\chi_r u), d(\chi_r u) \rangle - \langle d(\chi_r^2 u), du \rangle + \chi_r^2 u L_g u \right) dV \\
&= \int_M \left( u^2 |d\chi_r|^2 + \chi_r^2 u L_g u \right) dV \\
&= \int_{U_N(2r)} u^2 |d\chi_r|^2 dV + \int_M u L_g u \, dV + \int_{U_N(2r)} (\chi_r^2 - 1) u L_g u \, dV \\
&\leq \|u\|_{L^\infty(M)}^2 \frac{4}{r^2} \operatorname{vol}(U_N(2r)) + (\lambda - \varepsilon) \|u\|_{L^2(M)}^2 \\
&\quad + \|u\|_{L^\infty(M)} \|L_g u\|_{L^\infty(M)} \operatorname{vol}(U_N(2r)) \\
&\leq \left( \lambda - \varepsilon + C_1^2 \frac{4}{r^2} C_2 r^k + C_1^2 C_2 r^k \right) \|u\|_{L^2(M)}^2.
\end{aligned}$$

For the denominator of the Rayleigh quotient we have

$$\begin{aligned}
\|\chi_r u\|_{L^2(\widetilde{M})}^2 &= \|\chi_r u\|_{L^2(M)}^2 \\
&\geq \|u\|_{L^2(M)}^2 - \|u\|_{L^2(U_N(2r))}^2 \\
&\geq \|u\|_{L^2(M)}^2 - \|u\|_{L^\infty(M)}^2 \operatorname{vol}(U_N(2r)) \\
&\geq \left( 1 - C_1^2 C_2 r^k \right) \|u\|_{L^2(M)}^2.
\end{aligned}$$

This yields for the Rayleigh quotient

$$\frac{(L_{\widetilde{g}}(\chi_r u), \chi_r u)_{L^2(\widetilde{M})}}{\|\chi_r u\|_{L^2(\widetilde{M})}^2} \leq \frac{\lambda - \varepsilon + 4 C_1^2 C_2 r^{k-2} + C_1^2 C_2 r^k}{1 - C_1^2 C_2 r^k} \leq \lambda$$

for $r \leq R = R(\varepsilon, \Lambda, C_1^2, C_2, k)$ sufficiently small. From the unique continuation property of eigenfunctions it follows that the space $\{\chi_r u \mid u \in E_g(\lambda - \varepsilon)\}$ has the same dimension as $E_g(\lambda - \varepsilon)$ itself. Thus we have shown that the Rayleigh quotient of $L_{\widetilde{g}}$ is bounded by $\lambda$ on a space of dimension $\dim E_g(\lambda) = \mathcal{N}_g(\lambda)$. Hence

$$\mathcal{N}_{\widetilde{g}}(\lambda) \geq \mathcal{N}_g(\lambda - \varepsilon).$$

For the proof of this inequality the only assumption on $\widetilde{g}$ we have made is that $(\widetilde{M}, \widetilde{g})$ contains an isometric copy of $(M \setminus U_N(r), g)$ for sufficiently small $r$. We may therefore perform surgery inside $U_N(r)$ and choose an arbitrary extension $\widetilde{g}$ of $g|_{M \setminus U_N(r)}$ to the region in $\widetilde{M}$ replacing $U_N(r)$. In order to show

$$\mathcal{N}_{\widetilde{g}}(\lambda) \leq \mathcal{N}_g(\lambda + \varepsilon)$$

we will have to make more restrictive assumptions on $\widetilde{g}$.

Let $S_0$ be a lower bound of the scalar curvature of $(M, g)$. Choose a constant $S_1$ so large that

$$S_1 > \min\{0, \Lambda/c\} \quad \text{and} \quad \frac{\Lambda/c - S_0}{S_1 - \Lambda/c}\Lambda \leq \frac{\varepsilon}{2}. \tag{8}$$

By Proposition 2.1 of [1] there is a metric $g'$ on $M$ arbitrarily close to $g$ in the $C^1$-topology such that for $\mathrm{Scal} = \mathrm{Scal}_{g'}$

$$\mathrm{Scal} \geq \begin{cases} S_0 & \text{on all of } M, \\ 2S_1 & \text{on a neighborhood } U_0 \text{ of } N. \end{cases} \tag{9}$$

Since by Lemma 3.4 the eigenvalues of $L_g$ depend continuously on $g$ in the $C^1$-topology we may without loss of generality assume that (9) holds for $\mathrm{Scal} = \mathrm{Scal}_g$.

Next we choose $\eta > 0$ so small that

$$\frac{S_1 - S_0}{S_1 - \Lambda/c}((\Lambda + 1 - cS_0)\eta + \eta^2) \leq \frac{\varepsilon}{2}. \tag{10}$$

Now choose $r > 0$ so small that

- $2\sqrt{10}(\mathcal{N}_g(\Lambda + \varepsilon) + 1)r^{1/4} < \eta$,
- $U_N((2r)^{1/11}) \subset U_0$,
- $(2r)^{1/11}$ is no larger than the $R$ in Lemma 3.7.

We perform surgery along $N$ in the neighborhood $U_N(r)$. Hence $\widetilde{M}$ is of the form $\widetilde{M} = (M \setminus U_N(r)) \cup \widetilde{U}$. Surgery in codimension $\geq 3$ does not decrease scalar curvature too much if the metric $\widetilde{g}$ on $\widetilde{U}$ is chosen properly, see [10, Proof of Theorem A] and [16, Proof of Theorem 3.1]. We may assume

- $\mathrm{Scal}_{\widetilde{g}} \geq S_0$ on all of $\widetilde{M}$,
- $\mathrm{Scal}_{\widetilde{g}} \geq S_1$ on $\widetilde{U}$.

Now assume $\mathcal{N}_{\widetilde{g}}(\lambda) > m := \mathcal{N}_g(\lambda + \varepsilon)$ for some $\lambda \leq \Lambda$. This means there is an $(m+1)$-dimensional space $\widetilde{\mathcal{H}}$ of functions on $\widetilde{M}$ spanned by eigenfunctions of $L_{\widetilde{g}}$ for eigenvalues $\leq \lambda$. The space

$$\mathcal{H} := \{\chi_r u \mid u \in \widetilde{\mathcal{H}}\}$$

has the same dimension $m+1$ as $\widetilde{\mathcal{H}}$ by the unique-continuation property of eigenfunctions. We consider the elements of $\mathcal{H}$ as functions on $M$ (identically $0$ on $U_N(r)$). We will show that

$$\frac{(L_g v, v)_{L^2(M)}}{\|v\|^2_{L^2(M)}} \leq \lambda + \varepsilon \tag{11}$$

for all $v \in \mathcal{H}$. Hence $\mathcal{N}_g(\lambda + \varepsilon) \geq m + 1$, a contradiction. Let $v = \chi_r u \in \mathcal{H}$, $v \neq 0$. By Lemma 3.6 we have

$$\|u\|^2_{L^2(\widetilde{U} \cup A_N(r,2r))} \leq \frac{\Lambda/c - S_0}{S_1 - S_0}\|u\|_{L^2(\widetilde{M})}.$$

Thus

$$\begin{aligned}
\|v\|^2_{L^2(M)} &= \|\chi_r u\|^2_{L^2(M)} \geq \|u\|^2_{L^2(M\setminus U_N(2r))} \\
&\geq \left(1 - \frac{\Lambda/c - S_0}{S_1 - S_0}\right)\|u\|^2_{L^2(\widetilde{M})} \\
&= \frac{S_1 - \Lambda/c}{S_1 - S_0}\|u\|^2_{L^2(\widetilde{M})}
\end{aligned} \quad (12)$$

We write $u = u_1 + \cdots + u_{m+1}$ with $L_{\widetilde{g}} u_j = \mu_j u_j$, $\mu_j \leq \lambda$. Some of the $u_j$ may be zero. Fix $\rho \in [r, (2r)^{1/11}]$. Set $\widehat{M}_\rho := \widetilde{U} \cup A_N(r, \rho)$. Then $\widehat{M}_\rho \subset \widetilde{M}$ is a compact manifold with boundary $\partial \widehat{M}_\rho = S_N(\rho)$ and $\mathrm{Scal} \geq S_1$ on $\widehat{M}_\rho$. From

$$\begin{aligned}
\mu_j \|u_j\|^2_{L^2(\widehat{M}_\rho)} &= (L_g u_j, u_j)_{L^2(\widehat{M}_\rho)} \\
&= \int_{\widehat{M}_\rho} \langle \Delta u_j, u_j \rangle dV + c \int_{\widehat{M}_\rho} \mathrm{Scal}\, u_j^2 dV \\
&= \int_{\widehat{M}_\rho} \langle du_j, du_j \rangle dV - \int_{\partial \widehat{M}_\rho} u_j \partial_\nu u_j dA + c \int_{\widehat{M}_\rho} \mathrm{Scal}\, u_j^2 dV \\
&\geq -\int_{S_N(\rho)} u_j \partial_\nu u_j dA + cS_1 \|u_j\|^2_{L^2(\widehat{M}_\rho)}.
\end{aligned}$$

we get

$$\begin{aligned}
\int_{S_N(\rho)} u_j \partial_\gamma u_j dA &\geq (cS_1 - \mu_j)\|u_j\|^2_{L^2(\widetilde{M}_\rho)} \\
&\geq (cS_1 - \Lambda)\|u_j\|^2_{L^2(\widehat{M}_\rho)} \\
&\geq 0.
\end{aligned}$$

Lemma 3.7 yields

$$\begin{aligned}
\|u_j\|^2_{L^2(A_N(r,2r))} &\leq 10\, r^{5/2} \|u_j\|^2_{L^2(A_N(r,(2r)^{1/11}))} \\
&\leq 10\, r^{5/2} \|u_j\|^2_{L^2(\widetilde{M})} \\
&\leq 10\, r^{5/2} \|u\|^2_{L^2(\widetilde{M})}
\end{aligned}$$

and therefore

$$\begin{aligned}
\|u\|_{L^2(A_N(r,2r))} &= \|u_1 + \cdots + u_{m+1}\|_{L^2(A_N(r,2r))} \\
&\leq \sum_{j=1}^{m+1} \|u_j\|_{L^2(A_N(r,2r))} \\
&\leq \sqrt{10}(m+1)r^{5/4}\|u\|_{L^2(\widetilde{M})}.
\end{aligned}$$

Thus

$$\frac{2}{r}\|u\|_{L^2(A_N(r,2r))} \leq \eta \|u\|_{L^2(\widetilde{M})}. \quad (13)$$

Now we see

$$\begin{aligned}
(\Delta v, v)_{L^2(M)} &= \|d(\chi_r u)\|^2_{L^2(M)} \\
&= \|\chi_r du + u d\chi_r\|^2_{L^2(M)} \\
&\leq \left(\|\chi_r du\|_{L^2(M)} + \|u d\chi_r\|_{L^2(M)}\right)^2 \\
&\leq \left(\|du\|_{L^2(\widetilde{M})} + \frac{2}{r}\|u\|_{L^2(A_N(r,2r))}\right)^2 \\
&\stackrel{(13)}{\leq} \left(\|du\|_{L^2(\widetilde{M})} + \eta\|u\|_{L^2(\widetilde{M})}\right)^2 \\
&\leq (1+\eta)\|du\|^2_{L^2(\widetilde{M})} + (\eta + \eta^2)\|u\|^2_{L^2(\widetilde{M})}. \quad (14)
\end{aligned}$$

The scalar curvature of $\widetilde{M}$ satisfies $\mathrm{Scal} \geq S_1 \geq 0$ on $\widetilde{U} \cup U_N(2r)$ and therefore

$$(\mathrm{Scal}_g v, v)_{L^2(M)} = (\mathrm{Scal}_{\widetilde{g}} \chi_r u, \chi_r u)_{L^2(\widetilde{M})} \leq (\mathrm{Scal}_{\widetilde{g}} u, u)_{L^2(\widetilde{M})} \quad (15)$$

Combining (14) and (15) we obtain

$$\begin{aligned}
(L_g v, v)_{L^2(M)} &= (\Delta v, v)_{L^2(M)} + c(\mathrm{Scal}_g v, v)_{L^2(M)} \\
&\leq (1+\eta)\|du\|^2_{L^2(\widetilde{M})} + (\eta + \eta^2)\|u\|^2_{L^2(\widetilde{M})} + c(\mathrm{Scal}_{\widetilde{g}} u, u)_{L^2(\widetilde{M})} \\
&= (1+\eta)(L_{\widetilde{g}} u, u)_{L^2(\widetilde{M})} - \eta c(\mathrm{Scal}_{\widetilde{g}} u, u)_{L^2(\widetilde{M})} + (\eta + \eta^2)\|u\|^2_{L^2(\widetilde{M})} \\
&\leq (1+\eta)\lambda\|u\|^2_{L^2(\widetilde{M})} - \eta c S_0 \|u\|^2_{L^2(\widetilde{M})} + (\eta + \eta^2)\|u\|^2_{L^2(\widetilde{M})} \\
&\leq \left[\lambda + (\Lambda + 1 - cS_0)\eta + \eta^2\right] \|u\|^2_{L^2(\widetilde{M})}.
\end{aligned}$$

Dividing this by (12) yields

$$\begin{aligned}
\frac{(L_g v, v)_{L^2(M)}}{\|v\|^2_{L^2(M)}} &\leq \frac{S_1 - S_0}{S_1 - \Lambda/c}\left(\lambda + (\Lambda + 1 - cS_0)\eta + \eta^2\right) \\
&= \lambda + \frac{\Lambda/c - S_0}{S_1 - \Lambda/c}\lambda + \frac{S_1 - S_0}{S_1 - \Lambda/c}\left((\Lambda + 1 - cS_0)\eta + \eta^2\right) \\
&\leq \lambda + \frac{\varepsilon}{2} + \frac{\varepsilon}{2}
\end{aligned}$$

by (8) and (10). This proves (11) and we are done. □

## 4. Upper bound for $\kappa$

We are now going to use Corollary 3.2 together with results from bordism theory to prove an upper bound on $\kappa$ in terms of $\hat{A}$ for simply connected manifolds. The main technical tool is the deep result by Stolz [17] that a simply connected spin manifold of dimension at least five with vanishing $\alpha$-genus is spin bordant to a manifold with a positive scalar curvature metric. From Gromov and Lawson [10, Cor. C] we know that non-spin simply connected manifolds allow metrics of positive scalar curvature and hence have $\kappa = 0$.

We begin by showing that $\kappa$ is bounded by the number $p_n(\alpha)$ introduced in the following definition.

**Definition 4.1.** Let $n > 0$ be an integer and let $a \in \mathrm{KO}^{-n}(\mathrm{pt})$. If $a = 0$ let $p_n(a) = 0$, otherwise let $p_n(a)$ be the minimum of $|\pi_0(N)|$ for all compact, scalar flat spin manifolds $N$ with $\dim N = n$ and $\alpha(N) = a$.

Compact scalar flat manifolds of a given dimension and with a given value of the $\alpha$-genus can always be found, for example among manifolds with special holonomy, compare Theorems 4.4 and 4.6. So the minimum in the definition is a well-defined integer.

**Proposition 4.2.** *Let $M$ be a simply connected differentiable manifold of dimension $n \geq 5$. Then*
$$\kappa(M) \leq p_n(\alpha(M)).$$

*Proof.* If $M$ is non-spin or if $M$ is spin with $\alpha(M) = 0$ then $M$ admits a metric of positive scalar curvature and $\kappa(M) = 0$. Therefore let $M$ be spin and let $\alpha(M) \neq 0$. Let $N$ be a compact scalar flat spin manifold with $\dim N = \dim M$ and $\alpha(N) = \alpha(M)$. We then have $\alpha(M - N) = 0$ so by [17, Theorem B] there is a manifold $E$ such that $M$ is spin bordant to $E + N$ and $E$ has a metric of positive scalar curvature. The operator $L$ on $E + N$ has zero as an eigenvalue of multiplicity $|\pi_0(N)|$ and no negative eigenvalues. The $|\pi_0(N)| + 1^{\text{st}}$ eigenvalue is positive and can be assumed larger than 1.

Since $M$ is simply connected the spin bordism can be broken down into a sequence of elementary bordisms given by surgeries of codimension $\geq 3$, see [10, Proof of Thm. B]. From the Surgery Theorem 3.1 it follows that for each $\varepsilon > 0$ there is a metric on $M$ for which the first $|\pi_0(N)| + 1$ eigenvalues of $L$ on $M$ respectively $L$ on $E + N$ are $\varepsilon$-close. This means that the first $|\pi_0(N)|$ eigenvalues of $L$ on $M$ are have absolute value smaller than $\varepsilon$ while the $|\pi_0(N)| + 1^{\text{st}}$ one may still be assumed larger than 1. This implies $\kappa(M) \leq |\pi_0(N)|$ so we are done. □

**Remark 4.3.** The argument in the preceeding proof also shows that if $M$ and $M'$ are simply connected spin manifolds with $\dim M = \dim M' \geq 5$ and $\alpha(M) = \alpha(M')$ then $\kappa(M) = \kappa(M')$. So for simply connected spin manifolds $\kappa$ depends only on the dimension and on the $\alpha$-genus.

The next step is to find a concrete bound on $p_n(\alpha)$ in terms of $\alpha$. We need only consider manifolds $M$ of dimension equal to $0, 1, 2, 4 \mod 8$, since otherwise $\alpha(M) = 0$. We begin by looking at the case of dimension divisible by 4.

A scalar flat spin manifold with $\alpha \neq 0$ must have special holonomy, and its $\alpha$-genus is determined by the holonomy group. For the irreducible holonomy groups in dimension $n = 4m$ the cases of non-zero $\hat{A}$-genus are

- $\hat{A} = 2$ if $\text{Hol} = \text{SU}(2m)$,
- $\hat{A} = m + 1$ if $\text{Hol} = \text{Sp}(m)$,
- $\hat{A} = 1$ if $m = 2$ and $\text{Hol} = \text{Spin}(7)$,

and these cases all occur for compact manifolds, see for instance [13, Theorem 3.6.5]. The manifolds $N$ in Definition 4.1 can thus be taken as disjoint unions of products of manifolds with holonomy from this list, and $\hat{A}(N)$ is given by the sum of the products of the corresponding values of $\hat{A}$. The problem to determine $p_n(a)$ is to find a realisation of a given value of $\hat{A}$ with as few as possible terms in this sum.

We now give an upper bound on $p_n(a)$ using manifolds of dimensions 4 and 8 as generators. Let K3 be the K3-surface with a Ricci flat metric, this has $\hat{A}(\text{K3}) = 2$. Let $V_i$, $i = 0, \ldots, 4$ be spin 8-manifolds with Ricci flat metrics and $\hat{A}(V_i) = i$.

**Theorem 4.4.** *Let $M$ be a simply connected differentiable manifold of dimension $n \equiv 0$ mod $4$. Write $n = 8l$ or $n = 8l + 4$ with $l \geq 1$ and let $|\alpha(M)| = 4^l p + q$, $p \geq 0$, $0 \leq q < 4^l$. Then*
$$\kappa(M) \leq p + \min\{q, l\}.$$

This theorem follows from Proposition 4.2 and the following Lemma.

**Lemma 4.5.** *Let $n = 8l$ or $n = 8l + 4$, $l \geq 1$, and let $a$ be an integer. Write $|a| = 4^l p + q$, $0 \leq p$, $0 \leq q < 4^l$. Then*
$$p_n(a) \leq p + \min\{q, l\}.$$

*Proof.* For simplicity we assume $a \geq 0$. If $n = 8l$ we have
$$\alpha(pV_4^l + qV_1^l) = \hat{A}(pV_4^l + qV_1^l) = 4^l p + q = a$$
and
$$p_n(a) \leq |\pi_0(pV_4^l + qV_1^l)| = p + q.$$
If $n = 8l + 4$ we have
$$\alpha(\mathrm{K3} \times (pV_4^l + qV_1^l)) = \hat{A}(pV_4^l + qV_1^l) = 4^l p + q = a$$
so again
$$p_n(a) \leq |\pi_0(\mathrm{K3} \times (pV_4^l + qV_1^l))| = p + q.$$
To prove the inequality $p_n(a) \leq p + l$ set
$$W_i = V_4^i \times V_1^{l-1-i}, \quad i = 0, \ldots, l - 1.$$
Then $\dim W_i = 8l - 8$ and $\hat{A}(W_i) = 4^i$. Write $q = \sum_{i=0}^{l-1} q_i 4^i$, $0 \leq q_i < 4$. If $n = 8l$ we have
$$\alpha(pV_4^l + \sum_{i=0}^{l-1} V_{q_i} \times W_i) = \hat{A}(pV_4^l + \sum_{i=0}^{l-1} V_{q_i} \times W_i) = a$$
so
$$p_n(a) \leq |\pi_0(pV_4^l + \sum_{i=0}^{l-1} V_{q_i} \times W_i)| = p + l.$$
If $n = 8l + 4$ we have
$$\alpha(\mathrm{K3} \times (pV_4^l + \sum_{i=0}^{l-1} V_{q_i} \times W_i)) = \hat{A}(pV_4^l + \sum_{i=0}^{l-1} V_{q_i} \times W_i) = a$$
and
$$p_n(a) \leq |\pi_0(\mathrm{K3} \times (pV_4^l + \sum_{i=0}^{l-1} V_{q_i} \times W_i))| = p + l.$$
□

The inequality in Lemma 4.5 is in general not an equality. Therefore Theorem 4.4 can still be improved.

In dimensions $n \equiv 1, 2 \mod 8$ we have $\alpha(M) \in \mathrm{KO}^{-n}(\mathrm{pt}) \cong \mathbb{Z}/2\mathbb{Z}$. By $|\alpha(M)| \in \mathbb{Z}$ we mean $0$ if $\alpha(M)$ is zero and $1$ otherwise.

**Theorem 4.6.** *Let $M$ be a simply connected spin manifold of dimension $n = 8l + 1$ or $8l + 2$, $l \geq 1$. Then*
$$\kappa(M) = |\alpha(M)|.$$

*Proof.* If $\alpha(M) = 0$ then there is a metric of positive scalar curvature on $M$ and $\kappa(M) = 0$. If $\alpha(M) = 1$ then

$$\alpha\left(V_1^l \times (S^1)^a\right) = \alpha(M)$$

where $n = 8l + a$, $a = 1, 2$ and $S^1$ is the circle with the non-bounding spin structure. From Proposition 4.2 it follows that $\kappa(M) \leq 1$. Since $M$ does not allow a metric of positive scalar curvature we must have $\kappa(M) = 1$. □

## 5. MISCELLANEA

### 5.1. The stable limit.
Theorem 2.4 shows that $\kappa(M)$ can become arbitrarily large because this is true for the $\hat{A}$-genus. In a stable sense it takes only the values $0$ and $1$. More precisely, let $B$ be a compact simply connected 8-dimensional spin manifold with $\hat{A}(B) = 1$. Then $\alpha(M \times B) = \alpha(M)$ for all manifolds $M$. Recall that compact manifolds with holonomy $\mathrm{Spin}(7)$ are examples of such manifolds $B$.

**Theorem 5.1.** *Let $M$ be a simply connected spin manifold. Then*

$$\kappa(M \times B^p) \leq 1$$

*for all sufficiently large $p$.*

*Proof.* We have already seen that $\kappa \leq 1$ if the dimension is not divisible by four, so let $M$ be a compact simply connected spin manifold of dimension $n = 4m$. Since $\alpha(M \times B^p) = \alpha(M)$ we need to find a connected scalar flat manifold $N$ so that $\alpha(N) = \alpha(M)$ and $\dim N = 4m + 8p$ for some $p \geq 0$.

For $i \geq 1$ let $H_i$ be a compact $4i$-dimensional manifold with holonomy $Sp(i)$, so that $\hat{A}(H_i) = i + 1$. Assume $\hat{A}(M) \geq 0$.

First suppose $\dim(M) = 8l$ and set $\hat{A}(M) = 2^a(2b+1)$. If $a + b \leq l$ then

$$\hat{A}\left(V_1^{l-a-b} \times V_2^a \times H_{2b}\right) = \hat{A}(M)$$

so $\kappa(M) \leq 1$. If $a + b \geq l$ then

$$\hat{A}\left(V_2^a \times H_{2b}\right) = \hat{A}\left(M \times B^{a+b-l}\right)$$

and $\kappa(M \times B^{a+b-l}) \leq 1$.

Second suppose $\dim(M) = 8l + 4$ and set $\hat{A}(M) = 2 \cdot 2^a(2b+1)$. If $a + b \leq l$ then

$$\hat{A}\left(\mathrm{K3} \times V_1^{l-a-b} \times V_2^a \times H_{2b}\right) = \hat{A}(M)$$

so $\kappa(M) \leq 1$. If $a + b \geq l$ then

$$\hat{A}\left(\mathrm{K3} \times V_2^a \times H_{2b}\right) = \hat{A}\left(M \times B^{a+b-l}\right)$$

and $\kappa(M \times B^{a+b-l}) \leq 1$. □

## 5.2. $\kappa$ and coverings.
We now examine if there is any monotonicity property of $\kappa(M)$ when one replaces $M$ by a finite covering space (or conversely, by a finite quotient). It will turn out that this is not the case.

First we see that sometimes $\kappa(M)$ does not change when passing to a finite covering space. For example, if $M$ carries a metric of positive scalar curvature, then so does every covering space $\widetilde{M}$, i. e. $\kappa(M) = \kappa(\widetilde{M}) = 0$. Similarly, if $M = T^n$ is a torus, then for any $k$ positive and integral $M$ has $k$-fold coverings $\widetilde{M}$, again diffeomorphic to $T^n$. We then have $\kappa(M) = \kappa(\widetilde{M}) = 1$.

Next we see that $\kappa(M)$ can increase when we pass to a finite cover $\widetilde{M}$. For example, look at $M = \text{K3}\,\#T^4$. By Corollary 3.2 and (3) we have $\kappa(M) \leq \kappa(\text{K3}+T^4) = \kappa(\text{K3}) + \kappa(T^4) = 2$. Choose $k \geq 3 > \kappa(M)$. Any $k$-fold covering of $T^4$ yields a $k$-fold covering $\widetilde{M}$ of $M$ diffeomorphic to $\text{K3}\,\#\overset{k}{\cdots}\#\,\text{K3}\,\#T^4$. Then by Theorem 2.4

$$\kappa(\widetilde{M}) \geq 2^{-1} \cdot \hat{A}(\widetilde{M}) = 2^{-1} \cdot k \cdot \hat{A}(\text{K3}) = k > \kappa(M).$$

Finally, we observe that $\kappa(M)$ can also decrease when we pass to a covering. To construct examples pick an exotic sphere $\Sigma^n$ of dimension $n \equiv 1 \mod 8$ with non-trivial $\alpha(\Sigma^n)$. Such exotic spheres always exist in these dimensions. Let $k$ denote the order of $\Sigma^n$ in the (finite) group of $n$-dimensional manifolds homeomorphic to $S^n$. Put $M := (S^3/(\mathbb{Z}/k\mathbb{Z}) \times S^{n-3})\#\Sigma^n$. Then $\alpha(M) = \alpha(\Sigma^n)$ is non-trivial and since $M$ is spin it does not have a metric of positive scalar curvature. Thus $\kappa(M) \geq 1$. The universal covering $\widetilde{M}$ of $M$ is diffeomorphic to $(S^3 \times S^{n-3})\#\Sigma^n\#\overset{k}{\cdots}\#\Sigma^n = S^3 \times S^{n-3}$. Since this has a metric of positive scalar curvature we get $0 = \kappa(\widetilde{M}) < \kappa(M)$.

## 5.3. Non-positive eigenvalues.
For every Riemannian metric the operator $L$ has a finite number of non-positive eigenvalues. Similar to the invariant $\kappa(M)$ we can introduce the invariant $\kappa'(M)$ as the minimal number of non-positive eigenvalues of $L_g$ for all Riemannian metrics $g$ on the compact manifold $M$. Again this can be interpreted as a quantitative measure of how close a metric on $M$ can come to having positive scalar curvature. Trivially $\kappa' \leq \kappa$ and all results in the present work except Theorem 2.4 hold with $\kappa$ replaced by $\kappa'$. Unfortunately we do not know how to find lower bounds for $\kappa'$.

## 5.4. The operator $\Delta_g + c \cdot \text{Scal}_g$ for $0 < c < c_n$.
Let $c_n = \frac{n-2}{4(n-1)}$. The Conformal Laplacian $L_g = \Delta_g + c_n \cdot \text{Scal}_g$ on a manifold $M$ is positive if and only if $M$ admits a metric of positive scalar curvature, which is a topological condition on $M$. In [9, p. 45] the following question is posed.

**Question 5.2.** *What is the significance of the positivity of the operator $\Delta_g + c \cdot \text{Scal}_g$ for $0 < c < c_n$?*

(In the reference there is a mistake concerning the coefficient $c_n$ making the question slightly different and making the Observation on the same page irrelevant.)

We will now show that the condition in Question 5.2 gives no topological restrictions.

**Proposition 5.3.** *Let $M$ be a compact manifold of dimension $n \geq 3$ and let $0 < c < c_n$. Then there is a metric $\overline{g}$ on $M$ for which $\Delta_{\overline{g}} + c \cdot \text{Scal}_{\overline{g}} > 0$.*

*Proof.* Let $l$ be a function on $M$ whose critical points form a discrete set and let $g$ be a metric with $\mathrm{Scal}_g > 1$ on an open set containing the critical points of $l$. Set $\overline{g} = f^{\frac{4}{n-2}} g$ where $f = e^{tl}$ and $t$ is large enough so that

$$\frac{c}{c_n}\left(1 - \frac{c}{c_n}\right)|df|^2 + c \cdot \mathrm{Scal}_g \cdot f^2 = e^{2tl}\left(\frac{c}{c_n}\left(1 - \frac{c}{c_n}\right)t^2|dl|^2 + c \cdot \mathrm{Scal}_g\right)$$

is everywhere positive. From (1) and (2) it follows that

$$\begin{aligned}
&\int_M u\left(\Delta_{\overline{g}} + c \cdot \mathrm{Scal}_{\overline{g}}\right) u \, dV_{\overline{g}} \\
&= \int_M \left|f du + \frac{c}{c_n} u df\right|^2 + \left[\frac{c}{c_n}\left(1 - \frac{c}{c_n}\right)|df|^2 + c \cdot \mathrm{Scal}_g f^2\right] u^2 \, dV_g \\
&> 0
\end{aligned}$$

for all $u \in C^\infty(M)$, $u \not\equiv 0$, so $\Delta_{\overline{g}} + c \cdot \mathrm{Scal}_{\overline{g}} > 0$. $\square$

This also shows that although the surgery result in Theorem 3.1 holds for $\Delta + c \cdot \mathrm{Scal}$, $c > 0$, we would get a trivial invariant if we defined the kappa-invariant using this operator with $0 < c < c_n$.

## References

1. C. Bär and M. Dahl, *Surgery and the spectrum of the Dirac operator*, to appear in J. Reine Angew. Math.
2. H. Baum, *Spin-Strukturen und Dirac-Operatoren über pseudoriemannschen Mannigfaltigkeiten*, BSB B. G. Teubner Verlagsgesellschaft, Leipzig, 1981.
3. T. Branson, *Kato constants in Riemannian geometry*, Math. Res. Lett. **7** (2000), 245–261.
4. D. Calderbank, P. Gauduchon, and M. Herzlich, *Refined Kato inequalities and conformal weights in Riemannian geometry*, J. Funct. Anal. **173** (2000), 214–255.
5. J. Cheeger, *A lower bound for the smallest eigenvalue of the Laplacian*, Probl. Analysis, Sympos. in Honor of Salomon Bochner, Princeton Univ. 1969, 1970, pp. 195–199.
6. Y. Colin de Verdière, *Construction de laplaciens dont une partie finie du spectre est donnée*, Ann. Sci. Éc. Norm. Supér., IV. Sér. **20** (1987), 599–615.
7. A. Futaki, *Scalar-flat closed manifolds not admitting positive scalar curvature metrics*, Invent. Math. **112** (1993), 23–29.
8. S. Gallot and D. Meyer, *D'un résultat hilbertien à un principe de comparaison entre spectres. Applications*, Ann. Sci. École Norm. Sup. **21** (1988), 561–591.
9. M. Gromov, *Positive curvature, macroscopic dimension, spectral gaps and higher signatures*, Functional analysis on the eve of the 21st century, Vol. II, Birkhäuser, Boston, MA, 1996, pp. 1–213.
10. M. Gromov and H. B. Lawson, *The classification of simply connected manifolds of positive scalar curvature*, Ann. of Math. **111** (1980), 423–434.
11. ———, *Positive scalar curvature and the Dirac operator on complete Riemannian manifolds*, Publ. Math., Inst. Hautes Etud. Sci. **58** (1983), 295–408.
12. O. Hijazi, *A conformal lower bound for the smallest eigenvalue of the Dirac operator and Killing spinors*, Commun. Math. Phys. **104** (1986), 151–162.
13. D. D. Joyce, *Compact manifolds with special holonomy*, Oxford University Press, Oxford, 2000.
14. V. Mathai, *Nonnegative scalar curvature*, Ann. Global Anal. Geom. **10** (1992), no. 2, 103–123.
15. J. Petean, *The Yamabe invariant of simply connected manifolds*, J. Reine Angew. Math. **523** (2000), 225–231.
16. J. Rosenberg and S. Stolz, *Metrics of positive scalar curvature and connections with surgery*, Surveys on surgery theory, Vol. 2, Princeton Univ. Press, Princeton, NJ, 2001, pp. 353–386.
17. S. Stolz, *Simply connected manifolds of positive scalar curvature*, Ann. of Math. **136** (1992), 511–540.

CHRISTIAN BÄR, UNIVERSITÄT HAMBURG, FB MATHEMATIK, BUNDESSTR. 55, 20146 HAMBURG, GERMANY

MATTIAS DAHL, INSTITUTIONEN FÖR MATEMATIK, KUNGL TEKNISKA HÖGSKOLAN, 100 44 STOCKHOLM, SWEDEN

*E-mail address*: baer@math.uni-hamburg.de, dahl@math.kth.se